\documentclass[12pt]{article}
\usepackage{amssymb,url,curves}
\usepackage{amsmath}
\usepackage{graphicx}

\newtheorem{theorem}{Theorem}[section]
\newtheorem{prop}[theorem]{Proposition}
\newtheorem{cor}[theorem]{Corollary}
\newtheorem{lemma}[theorem]{Lemma}
\newtheorem{conjecture}[theorem]{Conjecture}

\newcommand{\qed}{\hfill$\Box$}

\newcommand{\id}{\mathrm{id}}

\newtheorem{unnumber}{}

\newtheorem{prepf}[unnumber]{Proof}
\newenvironment{pf}{\prepf\rm}{\endprepf}

\newcommand{\Aut}{\mathop{\mathrm{Aut}}}

\begin{document}

\title{Trees and cycles}
\author{Peter J. Cameron and Liam Stott\\School of Mathematics and Statistics\\
University of St Andrews}
\date{Autumn 2020}
\maketitle

\begin{quote}
\dots\ the cycle has taken us up through forests
\begin{flushright}
Robert M. Pirsig~\cite{pirsig}
\end{flushright}
\end{quote}

\begin{abstract}
Let $T$ be a tree on $n$ vertices. We can regard the edges of $T$ as
transpositions of the vertex set; their product (in any order) is a cyclic
permutation. All possible cyclic permutations arise (each exactly once) if
and only if the tree is a star. In this paper we find the number of realised
cycles, and obtain some results on the number of realisations of each cycle, 
for other trees. We also solve the inverse problem of the number of trees
which give rise to a given cycle. On the way, we meet some familiar number
sequences including the Euler and Fuss--Catalan numbers.
\end{abstract}

\section{Introduction}

Let $T$ be a tree on the vertex set $\{1,\ldots,n\}$, with edge set $E(T)$.
We regard the edge $e\in E(T)$ which joins vertices $i$ and $j$ as the
transposition $(i,j)$ in the symmetric group $S_n$. These transpositions
generate the symmetric group, and form a \emph{minimal generating set} for
$S_n$, in the sense that no proper subset is a generating set.

Whiston~\cite{whiston} showed that the largest size of a minimal generating
set for $S_n$ is $n-1$, and Cameron and Cara~\cite{cc} showed that, for $n\ge7$,
there are just two types of minimal generating sets of size $n-1$, both derived
from trees: one consists of the set $E(T)$ as above, while the other has the
form $\{\{s\}\cup\{(st)^\epsilon:t\in E(T)\setminus\{s\},\epsilon=\pm1\}$,
for $s\in E(T)$, where $T$ is an arbitrary tree. The diameter of the Cayley
graph of $S_n$ with generating set $E(T)$ has been investigated by
Kraft~\cite{kraft}, and other properties of these graphs by Konstantinova
and co-authors (see for example~\cite{gkksv}).

The product of the transpositions corresponding to all edges of $T$ is an
$n$-cycle \cite{denes}. (This follows easily by induction from the fact that, if $g\in S_n$
has the property that $i$ and $j$ lie in different cycles, then these two
cycles are fused into a single cycle in the product $g(i,j)$.)

According to Cayley's Theorem, the number of trees on $\{1,\ldots,n\}$ is
$n^{n-2}$, while the number of orderings of the edges of such a tree is
$(n-1)!$. On the other hand, the number of $n$-cycles in $S_n$ is $(n-1)!$,
and all cycles are conjugate, so each cycle can be realised in $n^{n-2}$
ways as the product of all the edges in a tree \cite{denes}.

One might guess that, given any tree, every cycle arises uniquely from an
ordering of its edges; but a little thought shows that this is not so. Indeed,
the only trees with this property are the stars. This raises several questions:
\begin{enumerate}
\item Given a tree $T$, how many different cycles arise from multiplying its
edges in arbitrary order?
\item What is the distribution of the numbers of occurrences of cycles from
a given tree?
\item The inverse problem: given a cycle $c$, how many trees give rise to
$c$ when the edges are multiplied together in some order?
\end{enumerate}

To illustrate, Table~\ref{t6} shows the frequencies of cycles arising from the
six non-isomorphic trees on $6$ vertices. (The entry $(x,y)$ means that $y$
cycles have frequency $x$.) This and other computations were performed using
\textsf{GAP}~\cite{gap}. The column labelled Diameter gives the diameter of
the Cayley graph of $S_6$ with connection set $E(T)$. These values were
calculated using the \textsf{GAP} package GRAPE~\cite{grape}.

\begin{table}[htbp]
\begin{center}

\begin{tabular}{|c|c|c|c|}
\hline
Tree & Cycles & Diameter & Cycle frequencies
\\
\hline
\setlength{\unitlength}{1mm}
\begin{picture}(15,10)
\multiput(0,5)(3,0){6}{\circle*{1}}
\put(0,5){\line(1,0){15}}
\end{picture} 
& 16 & 15 &
$(1,2),(4,4),(6,2),(9,4),(11,2),(16,2)$
\\
\hline
\setlength{\unitlength}{1mm}
\begin{picture}(15,10)
\multiput(1,5)(3,0){4}{\circle*{1}}
\multiput(14,1)(0,8){2}{\circle*{1}}
\put(1,5){\line(1,0){9}}
\put(10,5){\line(1,1){4}}
\put(10,5){\line(1,-1){4}}
\end{picture}
& 24 & 11 &
$(1,4),(3, 4),(4,4),(6,4),(7,4),(9,4)$
\\
\hline
\setlength{\unitlength}{1mm}
\begin{picture}(15,10)
\multiput(0,5)(3.5,0){5}{\circle*{1}}
\put(7,1.5){\circle*{1}}
\put(0,5){\line(1,0){14}}
\put(7,1.5){\line(0,1){3.5}}
\end{picture}
& 24 & 11 &
$(1,2),(2,4),(3,4),(4,4),(7,4),(8,4),(11,2)$
\\
\hline
\setlength{\unitlength}{1mm}
\begin{picture}(15,10)
\multiput(2,1)(0,8){2}{\circle*{1}}
\multiput(14,1)(0,8){2}{\circle*{1}}
\multiput(6,5)(4,0){2}{\circle*{1}}
\multiput(2,1)(8,4){2}{\line(1,1){4}}
\multiput(2,9)(8,-4){2}{\line(1,-1){4}}
\put(6,5){\line(1,0){4}}
\end{picture}
& 36 & 10 &
$(1,8),(3,16),(4,4),(6,8)$
\\
\hline
\setlength{\unitlength}{1mm}
\begin{picture}(15,10)
\multiput(1,5)(4,0){4}{\circle*{1}}
\multiput(9,1)(0,8){2}{\circle*{1}}
\put(1,5){\line(1,0){12}}
\put(9,1){\line(0,1){8}}
\end{picture}
& 48 & 9 & 
$(1,12),(2,12),(3,12),(4,12)$
\\
\hline
\setlength{\unitlength}{1mm}
\begin{picture}(15,10)
\multiput(3,2)(0,6){2}{\circle*{1}}
\put(2,5){\circle*{1}}
\put(7,5){\circle*{1}}
\multiput(11,3)(0,4){2}{\circle*{1}}
\put(3,2){\line(4,3){4}}
\put(2,5){\line(1,0){5}}
\put(3,8){\line(4,-3){4}}
\put(7,5){\line(2,-1){4}}
\put(7,5){\line(2,1){4}}
\end{picture}
& 120 & 7 &
$(1,120)$
\\
\hline
\end{tabular}
\end{center}
\caption{\label{t6}Cycles from $6$-vertex trees}
\end{table}

In this paper, we answer the first and third question, and give some 
information about the second. One feature of the results, is that a couple of
famous integer sequences, the Euler numbers and the Fuss--Catalan  numbers,
come up in the investigation. More specifically:
\begin{itemize}
    \item the number of occurrences of the most frequent cycle obtained from the $n$-edge path is the $n$th Euler number $E_n$;
    \item the number of trees which realise a given $(n+1)$-cycle is the $n$th Fuss-Catalan number.
\end{itemize}

Let $O(T)$ denote the set of all orderings of the edges of a tree $T$, and let $C(T)$ denote the set of cycles arising from orderings of the edges of $T$; for $c\in C(T)$, the \emph{multiplicity} of $c$ is the number of
orderings of the edges of $T$ for which the product of the transpositions
is~$c$.

\section{Cycles from a tree}

In this section we answer the first question stated in the introduction: how many distinct cycles arise from a given tree $T$? The answer is summarised in Theorem 2.3. Before proceeding it is necessary to discuss a geometric interpretation of the cycles which arise as products of the edges of $T$, since the proof of the theorem will rely in part on a lemma which emerges from this interpretation. Indeed, this lemma turns out to be foundational to much of the work done in this paper.

First notice that the set-up of our problem is naturally reached via a geometric, or dynamical, interpretation. We consider the edges of $T$ as transpositions on $\{1,\ldots,n\}$: the edge $e=\{i,j\}$ corresponds to the transposition $(i,j)$. Set-theoretically, the edge is just a pair of elements while the transposition is a function which sends the element $i$ to $j$ and vice versa, fixing everything else. This may seem like a bit of a leap, but this gap can be bridged by appealing to the geometry of $T$. If we imagine a person standing at the vertex $i$, they could cross the edge $e$ and arrive at the vertex $j$; similarly they can cross in the opposite direction to reach $i$ from $j$. Of course, if this person were at any other vertex $v$ they cannot cross $e$; any attempt to do so would fail and they would remain at $v$.

This interpretation can naturally be extended to sequences of edges (correspondingly, products of transpositions). To take a basic example, suppose we have $e_1=\{i,j\}$ and $e_2=\{j,k\}$; then their corresponding transpositions are $(i,j),(j,k)$ and so the sequence $e_1e_2$ corresponds to the $3$-cycle $(i,k,j)$. As before, imagine a person standing at the vertex $i$, moving along the tree with respect to the sequence $e_1e_2$. First they cross $e_1$, bringing them to the vertex $j$, and then $e_2$, arriving at $k$. Now at $k$ they move with respect to $e_1e_2$ again; first they attempt to cross $e_1$ but cannot since $e_1$ is not incident to $k$, so they next cross $e_2$ and arrive at $j$. Finally, starting at $j$ they cross $e_1$ to arrive at $i$ since they now cannot cross $e_2$. Starting at any other vertices they will remain fixed, since neither $e_1$ nor $e_2$ are incident to them.

So by appealing to the geometry of a tree $T$ we can arrive at the correspondence between edges and transpositions, since the transposition corresponding to an edge $e$ describes the part that $e$ plays in the dynamics of $T$. As such we can use the transpositions to formalise and generalise the idea described above. Given a tree $T$ and some sequence of its edges $s=e_1\ldots e_m$ we define the \emph{$k$th step of the traversal of $T$ from $i$ with respect to $s$} to be the unique path $p_i^k$ from the vertex $i(e_1\ldots e_m)^{k-1}$ to $i(e_1\ldots e_m)^k$ where the edges $e_j$ are identified with their corresponding transpositions. We refer to these vertices as those \emph{hit by} or \emph{landed on by} the traversal. We then define the \emph{traversal of $T$ from $i$ with respect to $s$} to be the concatenation of each step from $i$ in the order they appear; $p_i=p_i^1\ldots p_i^r$ where $r$ is the smallest number such that $i(e_1\ldots e_m)^r=i$. Notice that $p_i$ is necessarily a circuit and if $c_i$ is the cycle of $e_1\ldots e_m$ containing $i$ when written in disjoint cycle form then $p_i$ corresponds to $c_i$ and $c_i=(i,t(p_i^1),t(p_i^2),\ldots,t(p_i^{r-1}))$ which are the vertices hit by the traversal from $i$, where $t(p)$ denotes the terminal vertex of a path $p$. Finally we define the \emph{traversal of $T$ with respect to $s$} to be the set of circuits $p_i$ for each $i$ and thus the traversal $\{p_{i_1},\ldots p_{i_p}\}$ corresponds to the permutation $e_1\ldots e_m=c_{i_1}\ldots c_{i_p}$.

We are concerned with orderings of the edges $O(T)$. These are sequences as above which are maximal without replacement. As stated in the introduction, the products of transpositions corresponding to orderings are $n$-cycles and so the traversal of each must be a single circuit which hits every vertex. In fact, one can use this correspondence with traversals to obtain an alternative proof that any ordering obtains an $n$-cycle; it can be shown independently that a traversal with respect to a given ordering of a tree must be a single circuit which lands on every vertex. Related to this is the following lemma, which will prove useful for solving question (c) stated in the introduction.

\begin{lemma}
    Let $T$ be a tree with vertex set $\{1,\ldots,n\}$ and let $p=e_1\ldots e_m$ be the traversal of $T$ with respect to some ordering from $O(T)$. Then for each $e\in E(T)$ there are precisely two distinct numbers $j_1,j_2\in\{1,\ldots, m\}$ such that $e=e_{j_1}=e_{j_2}$.
\end{lemma}

\begin{pf}
    We proceed by induction on $n$. Considering the path of length $1$ with vertices $i_1,i_2$, there is only one edge $e$ and so only one ordering. The first step of the traversal from $i_1$ is $e$; the second step starts at $i_2$ and is also $e$, returning to $i_1$. Thus the only traversal for an ordering for this tree is $ee$, satisfying the claim.
    
    Now assume the inductive hypothesis and suppose we are given an ordering $\sigma\in O(T)$, we find its traversal $p$ step by step. Pick any vertex to start at and call it $i$; we denote by $e_1,\ldots,e_d$ the edges incident to $i$ and assume they are labelled so that they appear in this order in $\sigma$. Notice that each $e_j$ leads to a subtree of $T$ which we denote $T_j$ and the $T_j$ are pairwise disjoint. The first step of the traversal is the path between $i$ and $i\sigma$; the first edge incident to $i$ appearing in $\sigma$ is $e_1$ and thus $i\sigma$ must be a vertex of $T_1$. Let $e_{1,1}\ldots e_{1,k}$ be the part of the first step of the traversal which is in $T_1$ (this may be empty) and let $v_1$ be the vertex such that $e_1=\{i,v_1\}$ (we may have $v_1=i\sigma$).
    
    By the inductive hypothesis each $e\in E(T_1)$ appears precisely twice in any traversal of $T_1$ with respect to an ordering of its edges; removing edges in $E(T)\setminus E(T_1)$ from $\sigma$, which we denote $\sigma|_{T_1}$, gives such an ordering and indeed, if $v\in V(T_1)$ then $v\sigma=v\sigma|_{T_1}$ as long as $v\sigma\in V(T_1)$. So we obtain a traversal $p_{T_1}$ where each $e\in E(T_1)$ appears precisely twice, hits each vertex of $T_1$ and ends at $i\sigma$ (since that is where it started). Notice that $e_{1,1}\ldots e_{1,k}$ must be a suffix of $p_{T_1}$ as once $e_{1,1}$ is crossed the rest must immediately follow as the edges appear in that order in $\sigma$ and we must then land on $i\sigma$ where the traversal ends (we know this because they are a suffix of the first step of the traversal $p$). Since $e_{1,1}$ is incident to $v_1$ the edge preceding it in $p_{T_1}$ must be an edge $e_{v_1}$ also incident to $v_1$ and this edge must be the rightmost edge incident to $v_1$ left of $e_{1,1}$ in $\sigma|_{T_1}$ or the rightmost edge incident to $v_1$ in $\sigma|_{T_1}$. Therefore, in $\sigma$, it is either the rightmost edge of $T_1$ left of $e_1$ incident to $v_1$ or (if there are no edges incident to $v_1$ left of $e_1$) it is the rightmost edge incident to $v_1$. In either case the final step of $p_{T_1}$ differs from the corresponding step of $p$; after $e_{v_1}$, instead of $e_{1,1}$ the traversal $p$ crosses $e_1$ and immediately $e_2$ without landing on $i$. So $p$ and $p_{T_1}$ disagree on the path $e_{1,1}\ldots e_{1,k}$, but these edges are still included in $p$ precisely twice so far since they were counted once in $p_{T_1}$ where they do agree, and once immediately before $p_{T_1}$ started. Thus we redefine $p_{T_1}$ by removing the suffix $e_{1,1}\ldots e_{1,k}$ and attaching it as a prefix; now $p_{T_1}$ is the part of $p$ containing all edges of $T_1$ which appear in $p$ and each appear precisely twice. In the case where $v=i\sigma$ we have that $p_{T_1}$ and $p$ agree apart from landing on $v_1$ the second time, and so no changes need be made to $p_{T_1}$.
    
    Now after landing on the final vertex of $T_1$, $p$ lands on a vertex of $T_2$ and we have precisely the same situation as before. Thus this continues inductively until we reach $T_d$; in the final step of $p$ we start at a vertex in $T_d$ and cross $e_d$. Since $e_d$ is the rightmost edge incident to $i$ in $\sigma$, $p$ must land on $i$ and the traversal is complete. We now see that $p=e_1p_{T_1}e_1e_2\ldots e_{d-1}e_dp_{T_d}e_d$ where each edge of $T_j$ appears precisely twice in $p_{T_j}$, as required. \qed
\end{pf}

Using similar ideas to the above we can prove the following lemma which, while a simple idea, turns out to be crucial throughout the paper and in particular will be used to prove the next theorem.

\begin{lemma}
    Let $\sigma,\tau\in O(T)$. Then $\sigma,\tau$ give the same $n$-cycle $c\in C(T)$ if and only if they differ only by some sequence of commuting non-adjacent edges.
\end{lemma}

\begin{pf}
    Since non-adjacent edges correspond to disjoint transpositions, it is trivial that if the orderings differ only by a sequence of commuting non-adjacent edges then they are equal in the symmetric group. We show that if they differ by a sequence of commutes which involves at least one pair of adjacent edges then they must be distinct.
    
    Indeed, suppose $\sigma=e_1\ldots e_je_{j+1}\ldots e_{n-1}\in O(T)$ and consider $\sigma'=e_1\ldots e_{j+1}e_j\ldots e_{n-1}$ where $e_j,e_{j+1}$ are incident to a common vertex denoted $i$. Listing the edges incident to $i$ by the order they appear in each ordering, we have $e_{i_1},\ldots,e_{i_k},e_j,e_{j+1},e_{i_{k+3}},\ldots,e_{i_d}$ for $\sigma$ and $e_{i_1},\ldots,e_{i_k},e_{j+1},e_j$, $e_{i_{k+3}},\ldots,e_{i_d}$ for $\sigma'$. We denote the subtree that each $e_m$ leads to by $T_m$ and the traversal of $T$ with respect to $\sigma$ (resp. $\sigma'$) by $p$ (resp. $p'$). As illustrated in the proof of Lemma 2.1 above, $p$ and $p'$ hits each vertex of $T_m$ in some order then each vertex of $T_{m+1}$ according to the order of their respective listing of edges incident to $i$. The order in which a traversal of $T$ with respect to an ordering in $O(T)$ hits the vertices fully determines the cycle in $C(T)$ obtained from that ordering. Thus we can see that $p$ hits the vertices in $T_j$ then those in $T_{j+1}$ while $p'$ hits the vertices in $T_{j+1}$ then those in $T_j$ and hence they must correspond to distinct $n$-cycles. This argument can be applied inductively so we see that a sequence of commutes involving any number of pairs of adjacent edges will result in a necessarily different cycle. \qed
\end{pf}

\begin{theorem}
Let $T$ be a tree with vertex set $\{1,\ldots,n\}$, and suppose that the
vertex $i$ has valency $d_i$. Then the number of cycles which arise from
multiplying together the edges of $T$ is $\prod_{i=1}^nd_i!$.
\label{t:t2c}
\end{theorem}

\begin{pf}
    Given a tree $T$ with vertex set $\{1,\ldots,n\}$ we show that the tree $T'$ obtained by attaching a leaf (edge and vertex) to a vertex $i$ of $T$ with degree $d$ has the property that $|C(T')|=(d+1)|C(T)|$. The theorem then follows by induction on $n$.
    
    Let $e_1,\ldots,e_d$ be the vertices incident to $i$ and let $\sigma\in O(T)$. Then $\sigma$ has the form \[\sigma=g_1e_1g_2\ldots g_de_dg_{d+1}\] where $g_j$ is a product of edges from $E(T)\setminus\{e_1,\ldots,e_d\}$ such that $g_j,g_k$ have no edges in common for $j\neq k$ and each $e\in E(T)\setminus\{e_1,\ldots,e_d\}$ is used by $g_j$ for some $j$. We allow the possibility $g_j=\id$ for any $j$.
    
    Let $e_{d+1}$ be the edge attached to $i$ to obtain $T'$. Notice that given an ordering $\sigma'\in O(T')$ we obtain an ordering $\sigma\in O(T)$ by removing $e_{d+1}$; as such, $\sigma'$ can be obtained from $\sigma$ by inserting $e_{d+1}$ in between the appropriate edges of $\sigma$. Further, given $\sigma\in O(T)$ we can obtain an ordering $\sigma'\in O(T')$ by inserting $e_{d+1}$ in between any two consecutive edges of $\sigma$.
    
    So consider inserting $e_{d+1}$ somewhere into $\sigma=g_1e_1g_2\ldots g_de_dg_{d+1}$; for a fixed $j$, if we insert $e_{d+1}$ in any two places so that it is beside an edge of $g_j$ the two orderings from $O(T')$ we obtain will give the same cycle by Lemma 2.2, since they differ only by a sequence of commutes of $e_{d+1}$ with edges from $g_j$. On the other hand, if we obtain $\sigma_j,\sigma_k\in O(T')$ by inserting $e_{d+1}$ so it is beside an edge of $g_j$ and $g_k$ respectively, for $j< k$, the cycle given by each must be distinct since to reach $\sigma_k$ from $\sigma_j$ by a sequence of commutes would require commuting $e_{d+1}$ and $e_{j+1}$, which are adjacent in $T$. Thus for every cycle in $C(T)$ there are $d+1$ cycles in $C(T')$ (one for each $g_j$) as required. \qed
\end{pf}

\begin{cor}
\begin{enumerate}
\item The only trees which give rise to all possible cycles are the stars.
\item The trees which give rise to the smallest number of cycles (namely
$2^{n-2}$) are the paths.
\end{enumerate}
\end{cor}

\begin{pf}
The degrees $d_i$ sum to $2(n-1)$. Since $(d_i-1)!(d_j+1)!>d_i!d_j!$ if
$d_i<d_j$, we maximise the number of cycles by moving degrees to the extreme
values $(n-1,1,1,\ldots,1)$, and minimise the number by moving them close
to the mean values, $(2,2,\ldots,2,1,1)$. \qed
\end{pf}

\section{Distributions}

We now turn to the second question: for a given tree $T$, what are the multiplicities of the cycles which arise? We obtain a characterisation of trees with cycles of multiplicity $1$ and a counting formula for the number of such cycles. We also find the number of cycles with extremal multiplicities for two classes of tree; paths and forked paths.

Before detailing our results we first develop another alternative view of the problem at hand, this time order theoretic, which will prove useful for what we show and may have uses for future work.

As before, consider a tree $T$ with vertex set $\{1,\ldots,n\}$ and take an ordering of its edges $\sigma=e_1\ldots e_{n-1}\in O(T)$. Notice that each such ordering corresponds to a linear order on $E(T)$ in the sense of a partial order where any two elements are comparable; if $<$ is the linear order corresponding to $\sigma$ then $e_1<e_2<\ldots<e_{n-1}$. As we know, multiple different orderings can give the same cycle, so what is it that unifies these orderings?

Lemma 2.2 tells us that -- as far as the cycle obtained from an ordering is concerned -- the only relevant thing is, for each vertex, what order do the edges incident to it appear in the ordering. We can see this fact at play in the proofs from Section 2. Of course, each edge is incident to two vertices so the constraints on how we can manipulate an ordering without changing the cycle is more complex than this, but those constraints emerge from this basic principle.

It is this idea which motivates the following definition. First we require some notation and terminology. When working with a partially ordered set we will usually consider the partial order ($p$) and the set it is defined over ($P$) separately, regarding $p$ as a subset of $P\times P$. However, when working with a linear order $l$ it will often be convenient to regard it as a bijection $l:P\rightarrow [m]$ where $m=|P|$ and we will switch freely between these two conceptions. For a partial order $p$ over $P$ we denote by $L(p)$ the set of linear extensions of $p$; these are the linear orders over $P$ which contain $p$ as a subset. For a partial order $p$ on a set $P$ we define the \emph{inverse of $p$} to be the partial order $p^{-1}=\{(y,x)\ |\ (x,y)\in p\}$ on $P$.

If $p$ is a partial order on a set $P$ then we say $q\subseteq p$ is a \emph{suborder} of $p$ if it a partial order on $P$ and in this case we call $p$ a \emph{refinement} of $q$. For $Q\subseteq P$ we say $q$ is the \emph{induced suborder} of $p$ on $Q$ if $q=p\cap (Q\times Q)$. If $p$ is a partial order on $P$ and $q$ is a partial order on a subset $Q\subseteq P$ then for ease we denote by $p\cup q$ the smallest partial order on $P$ which contains both $p$ and $q$; in other words, we implicitly take the transitive closure of the set $p\cup q$.

Now, let $\sigma\in O(T)$ be an ordering. For each vertex $i$ we consider $\sigma|_i$; this is $\sigma$ with all edges apart from those incident to $i$ removed. Understanding $\sigma\subseteq E(T)\times E(T)$ as a linear order, $\sigma|_i\subseteq \sigma$ is a linear order on the edges incident to $i$; indeed, $\sigma|_i$ is the induced suborder of $\sigma$ on the set of edges incident to $i$ and we refer to it as the \emph{local suborder} of $\sigma$ at $i$. Thus we define the \emph{partial order with respect to $\sigma$}, denoted by $p_\sigma$, to be the transitive closure of \[\bigcup_{1\leq i\leq n}\sigma|_i\]

Notice that $p_\sigma$ is a subset of $\sigma$ and thus $\sigma\in L(p_\sigma)$. Indeed, the partial order with respect to $\sigma$ is defined to preserve the order in $\sigma$ which the edges incident to a vertex appear in, for each vertex. As such, any linear extension $\sigma'$ of $p_\sigma$ will have the edges incident to a given vertex appearing in the same order as in $\sigma$; in other words, $\sigma$ and $\sigma'$ differ only by commuting non-adjacent edges and thus they give the same cycle, by Lemma 2.2. Conversely, if $\sigma'$ is an ordering which gives the same cycle as $\sigma$ then, to reach one from the other, no two edges incident to a common vertex can commute and so any such edges must appear in the same order in each; hence $p_\sigma=p_{\sigma'}$ and $\sigma$ and $\sigma'$ are linear extensions of the same partial order. The conclusion of this discussion is summarised in the following theorem.

\begin{theorem}
    Let $T$ be a tree. There is a one-to-one correspondence between $C(T)$ and the set $P(T)=\{p_\sigma\ |\ \sigma\in O(T)\}$ such that for a given $\sigma\in O(T)$ the cycle $c$ arising from $\sigma$ corresponds to the partial order $p_\sigma$. In addition, for each $p\in P(T)$ its linear extensions $L(p)$ correspond precisely to the orderings in $O(T)$ giving the cycle which $p$ corresponds to.
\end{theorem}

This correspondence can be exploited to find a slicker alternative proof of Theorem 2.3. More importantly, it gives an entirely new method for finding the multiplicity of cycles for a given tree; rather than counting orderings which give a cycle directly, we can find the partial order which a cycle corresponds to (easily done given an ordering from which the cycle arises) and count the number of linear extensions. Having said that, counting linear extensions is no mean feat; it is a classical problem in order theory known to be difficult in general. But, of course, the partial orders arising from a given tree $T$ form a special class (depending on the class of tree) so there is hope for progress. As we will see this method is also useful for related questions; for example, it can be used to show that the cycle $(1,3,5,\ldots,6,4,2)$ is the unique most frequent cycle (besides its inverse) for the path with vertex set $\{1,\ldots,n\}$ labelled in ascending order from one end to the other.

For now, we focus on cycles with multiplicity $1$. As noted, the star gives rise to every cyclic permutation exactly once. But
for any other tree, not every cycle arises, and so some cycles occur with
multiplicity greater than $1$. 

\subsection{Cycles with multiplicity~$1$}

In this section we determine all the trees $T$ for which $C(T)$ contains
cycles with multiplicity~$1$.

\begin{figure}[htbp]
\begin{center}
\setlength{\unitlength}{1mm}
\begin{picture}(70,15)
\multiput(10,10)(10,0){6}{\circle*{1}}
\multiput(0,6)(0,2){5}{\circle*{1}}
\multiput(70,6)(0,4){3}{\circle*{1}}
\multiput(18,0)(4,0){2}{\circle*{1}}
\put(30,0){\circle*{1}}
\multiput(46,0)(2,0){5}{\circle*{1}}
\put(0,10){\line(1,0){70}}
\put(0,6){\line(5,2){10}}
\put(0,8){\line(5,1){10}}
\put(0,12){\line(5,-1){10}}
\put(0,14){\line(5,-2){10}}
\put(18,0){\line(1,5){2}}
\put(22,0){\line(-1,5){2}}
\put(30,0){\line(0,1){10}}
\put(46,0){\line(2,5){4}}
\put(48,0){\line(1,5){2}}
\put(50,0){\line(0,1){10}}
\put(52,0){\line(-1,5){2}}
\put(54,0){\line(-2,5){4}}
\put(60,10){\line(5,-2){10}}
\put(60,10){\line(5,2){10}}
\end{picture}
\end{center}
\caption{\label{f:cat}A caterpillar}
\end{figure}
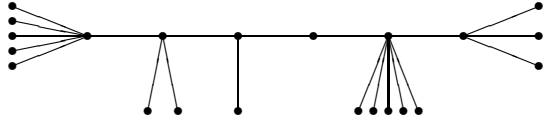

A \emph{caterpillar} is a tree with the property that removal of all leaves
gives rise to a path. This path is called the \emph{body} of the caterpillar.
Figure~\ref{f:cat} shows a caterpillar.

\begin{theorem}
A tree $T$ has the property that some cycle in $C(T)$ has multiplicity~$1$
if and only if $T$ is a caterpillar.
\end{theorem}

\begin{pf}
Suppose first that there is an ordering of the edges of $T$ such that the
product occurs with multiplicity~$1$ in $C(T)$. Then edges which are adjacent
in the ordering must meet at a vertex, else the transpositions would commute
and could be swapped. In other words, the ordering is a Hamiltonian path in
the line graph of $T$.

We claim that the tree with three paths of length~$2$ from a vertex is a
forbidden subgraph. For consider the tree with edges $\{v_1,v_2\}$,
$\{v_2,v_3\}$, $\{v_3,v_4\}$, $\{v_4,v_5\}$, $\{v_3,v_6\}$, $\{v_6,v_7\}$.
At some point in the sequence, without loss, we arrive at edge $\{v_2,v_3\}$
from $\{v_1,v_2\}$. If we visit $\{v_3,v_4\}$ and $\{v_3,v_6\}$ and leave
along $\{v_6,v_7\}$, we can never revisit $\{v_4,v_5\}$. But if $\{v_4,v_5\}$
follows $\{v_3,v_4\}$, then we can never return to $\{v_3,v_6\}$. (This is
the smallest tree which is not a caterpillar.)

We conclude that, at any vertex of $T$, all but (at most) two of the incident
edges are leaves. Thus, removing the leaves gives a path. So $T$ is a
caterpillar.

Conversely, given a caterpillar, order the edges so that, if $\{v_i,v_j\}$ and
$\{v_j,v_k\}$ are edges at vertex $v_j$ which are not leaves, we visit all the
leaves at $v_j$ between $\{v_i,v_j\}$ and $\{v_j,v_k\}$ in the sequence of
edges. \qed
\end{pf}

\begin{cor}
Let $T$ be a caterpillar whose body is the path $v_1,\ldots,v_r$. For 
$i=1,\ldots,r$, let $l_i$ be the number of leaves incident with $v_i$. Then
the number of cycles with multiplicity~$1$ in $C(T)$ is
\[
\begin{cases}
    2\prod_{i=1}^rl_i! & \text{if } r>1 \\
    l_1! & \text{if } r=1.
\end{cases}
\]
\label{c:uniq_cyc}
\end{cor}

\begin{pf}
The sequence of edges must use the edges of the body in turn. If it starts
at the $v_1$ end, it must traverse the leaves at $v_1$ in some order before
using $\{v_1,v_2\}$, then the leaves at $v_2$ in some order before using
$\{v_2,v_3\}$, and so on. If $r>1$ we double the number since we may start
at either end. \qed
\end{pf}

Every tree on $6$ vertices is a caterpillar, and Table~\ref{t6} agrees with
the above Corollary. As noted, the smallest tree which is not a caterpillar has
$7$ vertices; the frequency distribution of cycles for this tree is
$(3,12)$, $(9,12)$, $(15,12)$, $(33,12)$.

\subsection{Paths}

For stars, every cycle is realised with multiplicity~$1$.
The obvious next case to look at is the opposite extreme, the paths, which
realise the smallest number of cycles.

\subsubsection{Best and worst cycles}

Consider the $n$-vertex path, with edges numbered consecutively from $1$ to
$(n-1)$. The edge transpositions are the \emph{Moore--Coxeter generators}
of the symmetric group. Empirically we found that the most frequent cycles are
$(1,3,5,\ldots,6,4,2)$ and its inverse, and that the frequencies of these
cycles are the \emph{Euler numbers} $E_n$ (sequence A000111 in the
On-line Encyclopedia of Integer Sequences~\cite{oeis}),
having generating function $\sec(x)+\tan(x)$~\cite[p.149]{stanley}.
The sequence begins
\[1, 1, 1, 2, 5, 16, 61, 272, 1385, 7936, 50521, \dots\]

We prove that the number of realisations of this cycle is the Euler number in Theorem 3.5 and show that no other cycle apart from its inverse does better in Theorem 3.11.

The second most frequent cycles give the sequence
\[0, 0, 0, 1, 3, 11, 40, 181, 917, 5263, 33486, \dots\]
where we have put $0$ if only the maximum frequency occurs. This sequence is
not in the OEIS. Can we find further information about it?

\subsubsection{Realising the Euler numbers}

Rather than a direct proof of the formula, we give a characterisation of those
permutations of the edges of the path which give rise to the cycle
$(1,3,5,\ldots,6,4,2)$; these are identified with the inverses of the
permutations counted by the Euler numbers.

To simplify notation, we use instead the path with $n+1$ vertices, so that
we are looking at permutations on $\{1,\ldots,n\}$.

We need to distinguish between the active and passive forms of a permutation.
The \emph{passive form} of a permutation $\sigma$ on $\{1,\ldots,n\}$ is the
$n$-tuple $[a_1,a_2,\ldots,a_n]$ containing each element of $\{1,\ldots,n\}$
just once; the \emph{active form} is the map $\sigma$ which takes $i$ to
$a_i$ for $i=1,\ldots,n$. We usually write the active form in disjoint cycle
notation; we use square brackets for the passive form to avoid confusion with
the active form of a cyclic permutation. The \emph{inverse} of a permutation 
$\sigma$ is the permutation whose active form is the inverse map of the active
form of $\sigma$.

Multiplying a permutation on the left by a transposition $(i,j)$ has the
effect of interchanging the elements in positions $i$ and $j$ in its passive
form.

An \emph{up-down permutation} is one whose passive form satisfies
\[a_1<a_2>a_3<a_4>\cdots.\]
The \emph{canonical up-down permutation} is the permutation whose passive form
is $[1,n,2,n-1,3,\ldots]$. Its inverse has passive form $[1,3,5,\ldots,6,4,2]$.
The number of up-down permutations on $\{1,\ldots,n\}$ is the
Euler number \cite[p.149]{stanley}. These and the analogous down-up
permutations are also known as \emph{alternating permutations} or \emph{zig-zag
permutations}.

A transposition $(i,j)$ is \emph{acceptable} for a permutation
$[a_1,\ldots,a_n]$ if $|i-j|>1$ but $|a_i-a_j|=1$.

\begin{lemma}
From any up-down permutation $\sigma$, we can reach the canonical up-down
permutation by multiplying on the left by a sequence $t_m,\ldots,t_1$ of
transpositions, such that for each $i$, $t_{i-1}\cdots t_1\sigma$ is up-down
and $t_i$ is acceptable for this permutation.
\end{lemma}

\begin{pf}
Suppose first that $n$ is even.
Let $\sigma$ be an up-down permutation with passive form $[a_1,\ldots,a_n]$.
We call the even-numbered positions \emph{peaks} and the others \emph{troughs}.
A peak $i$ is \emph{low} if $a_i\le n/2$, while a trough $j$ is \emph{high} if
$a_j\ge 1+n/2$. Note that the numbers of low peaks and high troughs are
equal, and that a low peak and a high trough cannot be adjacent.

The first step is to reduce the number of low peaks (assuming it is nonzero).
Let $i$ be the low peak for which $a_i$ is maximal, and suppose that
$a_i+1=a_j$. Then $j$ is either a trough or a high peak. If it is a trough,
then the transposition $(i,j)$ swaps the entries in these positions; so
in the new permutation, $i$ is a peak (and is either high, or low but higher
than before) while $j$ is a trough. If $j$ is a high peak, then $a_i=n/2$. 
So we can raise the height of the peak at $i$ to $n/2$. Similarly, we can
lower the height of the lowest high trough to $n/2+1$. Then we can swap these
two to reduce the number of low peaks by one.

For example, consider $[1,3,2,5,4,8,6,7]$, where the second position is a low
peak and the seventh is a high trough. Applying $(2,5)$ takes us to
$[1,4,2,5,3,8,6,7]$; then applying $(4,7)$ gives $[1,4,2,6,3,8,5,7]$. Then
we may apply $(2,7)$ to give $[1,5,2,6,3,8,4,7]$, with no low peaks or high
troughs.

After finitely many steps of this type, we reach the case where there are
no low peaks or high troughs. Since the peaks are pairwise non-adjacent, we
can arrange them in descending order by a number of swaps, realised by left
multiplication by acceptable transpositions (using \emph{bubblesort});
similarly we can arrange the troughs in ascending order. The result is the
canonical up-down permutation. In the above example, the troughs are already
sorted, and we sort the peaks by applying in turn $(4,8)$, $(2,8)$, $(2,4)$,
$(2,6)$, and $(4,6)$.

If $n$ is odd, the number of troughs is one more than the number of peaks,
since position $n$ is a trough. We split into two cases. If the value of $i$
for which $a_i=(n+1)/2$ is a trough, we regard it as a low trough, and proceed
as before. If it is a peak, there must be a high trough; we can reduce the
height of the lowest high trough $j$ to $(n+3)/2$ by the above method. After
swapping these, $a_j=(n+1)/2$ and $j$ is a trough. \qed
\end{pf}

\begin{theorem}
Number the edges of the $(n+1)$-vertex path from $1$ to $n$. Then a permutation
$\sigma$ of the edges gives rise to the cycle $(1,3,5,\ldots,6,4,2)$ if and
only if $\sigma$ is the inverse of an up-down permutation.
\end{theorem}

\begin{pf}
A simple calculation shows that the inverse of the canonical up-down 
permutation gives rise to this cycle, which we call the \emph{canonical} cycle. 

We observe that multiplying $\sigma=[a_1,\ldots,a_n]$ on the left by the
transposition $(i,j)$ is the same as multiplying it on the right by the
conjugate $(i,j)^\sigma=(a_i,a_j)$; this satisfies the dual conditions
$|i-j|=1$ but $|a_i-a_j|>1$. The effect on the inverse permutation is to
multiply on the left by this transposition.

The corresponding edges of the tree are disjoint but are adjacent in the order
of the product; so they commute, and can be swapped. Thus such a multiplication
does not change the cycle given by the permutation.

Now, by the lemma, an arbitrary inverse
up-down permutation can be obtained from this one by post-multiplying by a 
sequence of such transpositions. This operation does not change the cycle.
We conclude that all inverses of up-down permutations give rise to the same cycle.

For the converse, we know that the numbers of up-down permutations and
of permutations realising the canonical cycle are equal, and since one set
is contained in the other, the sets are equal. \qed
\end{pf}

Note that the inverse of the canonical cycle is realised by the
similarly-defined \emph{down-up permutations}.

\subsubsection{The Euler numbers are maximal}

We have given a characterisation of the orderings of the edges of a path which give the canonical cycle $(1,3,5,\ldots,6,4,2)$ and using this we conclude the multiplicity of this cycle is the Euler number $E_n$ where $n+1$ is the number of vertices. Now we use Theorem 3.1 to show that no other cycle has a multiplicity higher than this, and this cycle and its inverse are the only ones to attain it.

First, since we are specifically concerned with paths, we can make this correspondence more precise as demonstrated in the following proposition. As before, we consider the path $P_n$ on $n+1$ vertices with edges consecutively labelled $1,\ldots,n$ from left to right.

\begin{prop}
    The set $\mathbf{P}_n=P(P_n)$ is precisely the set of partial orders on $n$ points whose Hasse diagrams are paths.
\end{prop}

Before proceeding with the proof we note that a Hasse diagram which is a path with $n$ vertices may be defined by taking a simple path with $n$ vertices and specifying an appropriate orientation on each of its edges; say $1$ and $-1$ for indicating an upward and downward slope respectively. Further, specifying any such orientation on the path defines such a Hasse diagram.

\begin{pf}
    Let $p\in\mathbf{P}_n$, then there is some $\sigma\in O(P_n)$ such that $p=p_\sigma$. To determine the Hasse diagram of $p$ we consider $\sigma|_v$ for each vertex $v$ of $P_n$; recall the edges of a Hasse diagram are the cover relations for its partial order so this is all we need consider. There are two leaves in $P_n$ (which the edges $1$ and $n$ are incident to); the remaining vertices $v_i$ have two edges incident to each, $i$ and $i+1$ for $1\leq i\leq n-1$. Thus there is an edge between the vertices of the Hasse diagram of $p$ labelled $i$ and $i+1$ (the direction of the edge depending on $\sigma|_{v_i}$) for each $i$ and no other edges. In other words, the Hasse diagram of $p$ is a path.
    
    On the other hand, we know there is a bijection between $C(P_n)$ and $\mathbf{P_n}$ and by Corollary 2.4(b) we have $|C(P_n)|=2^{n-1}$. Thus it suffices to show that the number of Hasse diagrams who are paths with $n$ vertices is also $2^{n-1}$. Indeed, by the comment before the proof we can count these Hasse diagrams by counting the possible orientations on the edges of a simple path with $n$ vertices. This path has $n-1$ edges and so there are $2^{n-1}$ different orientations. \qed
\end{pf}

We will thus refer to the elements of $\mathbf{P}_n$ as partial orders of \emph{path-type}. We will frequently identify these partial orders with their Hasse diagrams (as with partial orders in general) and in light of the above comment it is coherent, and will be convenient, to identify these with $(n-1)$-tuples whose entries are $\pm 1$ corresponding to each Hasse diagram's orientation.

Our aim is to show that the canonical cycle $(1,3,5,\ldots,6,4,2)$ attains the highest multiplicity among the cycles in $C(P_n)$, and only its inverse does at least as well. To do this we should find the partial order in $\mathbf{P}_n$ corresponding to the canonical cycle; for this, all we need is an ordering from $O(P_n)$ which gives this cycle. By Theorem 3.5 any inverse of an up-down permutation on $E(P_n)$ corresponds to such an ordering and in particular we may use the inverse of the canonical up-down permutation $\Sigma$, which has passive form $[1,3,5,\ldots,6,4,2]$. Since all odd numbers appear before any even numbers, it is clear that $\Sigma|_{v_i}=i(i+1)$ when $i$ is odd and $\Sigma|_{v_i}=(i+1)i$ when $i$ is even, giving $(i,i+1)$ and $(i+1,i)$ as cover relations for $p_\Sigma$ respectively. Thus the canonical cycle corresponds to $(1,-1,1,\ldots,(-1)^n)\in\mathbf{P}_n$; this is known as the $\emph{zig-zag partial order}$ and we denote it by $Z_n$ (Figure~\ref{f:zig}).

\begin{figure}[htbp]
    \begin{center}
        \setlength{\unitlength}{1mm}
        \begin{picture}(60,15)
            \multiput(5,2)(4,8){2}{\circle*{1}}
            \multiput(13,2)(4,8){2}{\circle*{1}}
            \put(21,2){\circle*{1}}
            \put(5,2){\line(1,2){4}}
            \put(9,10){\line(1,-2){4}}
            \put(13,2){\line(1,2){4}}
            \put(17,10){\line(1,-2){4}}
            \put(21,2){\line(1,2){2}}
            \put(25,6){\dots}
            \put(32,6){\line(1,-2){2}}
            \multiput(38,10)(8,0){2}{\circle*{1}}
            \multiput(34,2)(8,0){3}{\circle*{1}}
            \put(34,2){\line(1,2){4}}
            \put(38,10){\line(1,-2){4}}
            \put(42,2){\line(1,2){4}}
            \put(46,10){\line(1,-2){4}}
        \end{picture}
    \end{center}
    \caption{\label{f:zig}The Hasse diagram of a zig-zag partial order $Z_n$ for odd $n$}
\end{figure}

We now define an operation on $\mathbf{P}_n$ which we will then show necessarily increases the number of linear extensions of the partial order it is applied to, under the appropriate conditions. Then, given any path-type partial order, we show that there is some sequence of applying this operation under such conditions which brings us to the zig-zag partial order, whence we conclude no other path-type partial order does better.

\begin{figure}[htbp]
    \begin{center}
        \setlength{\unitlength}{1mm}
        \begin{picture}(50,22)
            \multiput(2,2)(3,6){4}{\circle*{1}}
            \put(14,14){\circle*{1}}
            \put(2,2){\line(1,2){9}}
            \put(14,14){\line(-1,2){3}}
            \put(18,11){\vector(1,0){8}}
            \put(20,13){$\Phi_3$}
            \multiput(30,2)(3,6){3}{\circle*{1}}
            \put(42,14){\circle*{1}}
            \put(39,8){\circle*{1}}
            \put(30,2){\line(1,2){6}}
            \put(36,14){\line(1,-2){3}}
            \put(39,8){\line(1,2){3}}
        \end{picture}
        
        \begin{picture}(80,20)
            \multiput(2,2)(3,6){3}{\circle*{1}}
            \multiput(8,14)(3,-6){3}{\circle*{1}}
            \put(2,2){\line(1,2){6}}
            \put(8,14){\line(1,-2){6}}
            \put(18,8){\vector(1,0){8}}
            \put(20,10){$\Phi_4$}
            \multiput(30,2)(3,6){3}{\circle*{1}}
            \multiput(39,8)(3,6){2}{\circle*{1}}
            \put(30,2){\line(1,2){6}}
            \put(36,14){\line(1,-2){3}}
            \put(39,8){\line(1,2){3}}
            \put(46,8){\vector(1,0){8}}
            \put(48,10){$\Phi_2$}
            \multiput(58,2)(3,6){2}{\circle*{1}}
            \multiput(64,2)(3,6){2}{\circle*{1}}
            \put(70,2){\circle*{1}}
            \put(58,2){\line(1,2){3}}
            \put(61,8){\line(1,-2){3}}
            \put(64,2){\line(1,2){3}}
            \put(67,8){\line(1,-2){3}}
        \end{picture}
    \end{center}
    \caption{\label{f:ex}Some right-side inversions}
\end{figure}
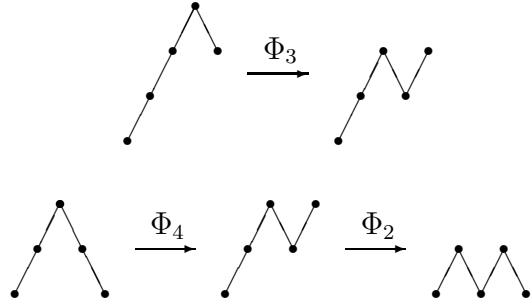

We define the operation $\Phi_i:\mathbf{P}_n\rightarrow\mathbf{P}_n$ for $1\leq i \leq n-1$ by \[(a_1,\ldots, a_{n-1})\Phi_i=(a_1,\ldots,a_{i-1},-a_i,-a_{i+1},\ldots,-a_{n-1})\] We call this operation \emph{right-side inversion at $i$}. Note that for $i=1$ the operation inverts $p\in\mathbf{P}_n$. A few examples of this operation being applied are shown in Figure~\ref{f:ex}; they begin with $(1,1,1,-1),(1,1,-1,-1)\in\mathbf{P}_5$.

The following two lemmas cover some useful results from order theory and are probably well known.

\begin{lemma}
    Let $p$ be a partial order on a set $P$. There is a bijection $\omega_p:L(p)\rightarrow L(p^{-1})$ and in particular $|L(p)|=|L(p^{-1})|$.
\end{lemma}

\begin{pf}
    We define $\omega_p:L(p)\rightarrow L(p^{-1})$ by $p'\omega_p=(p')^{-1}$. To see that $(p')^{-1}\in L(p^{-1})$ note that $(p')^{-1}$ contains $p^{-1}$ since $p'$ contains $p$. Since $(r^{-1})^{-1}=r$ for any partial order $r$ on $P$ it follows by symmetry that $\omega_p$ is invertible and hence a bijection. \qed
\end{pf}

\begin{lemma}
    Let $p$ be a partial order on some set $P$ and let $q$ be an induced suborder on a subset $Q\subseteq P$. Then $L(p)=\bigsqcup_{q'\in L(q)}L(p\cup q')$.
\end{lemma}

\begin{pf}
    The partial order $p\cup q'$ is a refinement of $p$ and so each $r\in L(p\cup q')$ is a linear extension of $p$ for all $q'\in L(q)$.
    
    Let $p'\in L(p)$. Take the induced suborder of $p'$ on $Q$ and denote it as $q'$ - notice that $q\subseteq q'$ since $q\subseteq p\subseteq p'$. Since $p'$ is a linear order any induced suborder is as well, so $q'\in L(q)$ and $p'$ is a linear order containing both $p$ and $q'$; thus $p'\in L(p\cup q')$.
    
    To see that the union is disjoint note that for $q',q''\in L(q)$ distinct there must exist $x,y\in Q$ so that $q'(x)<q'(y)$ and $q''(y)<q''(x)$ and consequently these relations must hold for any $r\in L(p\cup q')$ and $s\in L(p\cup q'')$. \qed
\end{pf}

The next lemma we require involves the \emph{Entringer numbers} $E(n,i)$ which are the number of down-up permutations on $\{1,\ldots,n+1\}$ whose first term is $i+1$ when written in passive form. They satisfy the recursion $E(n+1,i+1)=E(n+1,i)+E(n,n-i)$ for $0\leq i\leq n$ where $E(0,0)=1$ and $E(n,0)=0$ for all $n\geq 1$ \cite{stanley}.

\begin{lemma}
    Define $L(Z_n^{-1},i)=\{ z\in L(Z_n^{-1})\ |\ z(1)=i+1\}$. Then $|L(Z_n^{-1},i)|=E(n-1,i)$ for $0\leq i\leq n-1$.
\end{lemma}

\begin{pf}
    Recall that a down-up permutation on $\{1,\ldots,n\}$ is a permutation whose passive form $[a_1,\ldots,a_n]$ satisfies $a_1>a_2<a_3>\ldots$ and notice that a linear order $l$ on the set $\{1,\ldots n\}$ may be regarded as a permutation where the image of $l$ specifies how $l$ orders the elements of the set. In particular any $z\in L(Z_n^{-1})$ must satisfy $z(1)>z(2)<z(3)>\ldots$ and thus the number of ways of choosing the $z(j)$ such that $z(1)=i+1$ is the Entringer number $E(n-1,i)$. \qed
\end{pf}

\begin{prop}
    Let $p=(a_1,\ldots,a_{n-1})\in\mathbf{P}_n$ and let $i\geq 2$ be maximal such that $a_i=a_{i-1}$. Then $|L(p)|<|L(p\Phi_i)|$.
\end{prop}

\begin{pf}
    We assume without loss of generality that $a_i=a_{i-1}=-1$. Since inverting a partial order doesn't change its number of linear extensions, if $a_i=a_{i-1}=1$ we can invert $p$, apply right-side inversion at $i$ and then invert the resulting partial order.
    
    Consider the induced suborder $q$ of $p$ on $Q=\{i,i+1,\ldots,n\}$. Note that $q\cong Z_{n-i+1}^{-1}$ by maximality of $i$. Given $q'\in L(q)$ we can then consider the induced suborder $r_{q'}$ of $p\cup q'$ on $R=\{i-1,i,\ldots,n\}$.
    
    Similarly consider the induced suborder of $p\Phi_i$ on $Q$ (this is the inverse of $q$ so we denote it as $q^{-1}$) and further consider the induced suborder $r_{(q^{-1})'}$ of $p\Phi_i\cup (q^{-1})'$ on $R$ for $(q^{-1})'\in L(q^{-1})$. By Lemma 3.7 there is a bijection $\omega_q:L(q)\rightarrow L(q^{-1})$ which sends a linear extension of $q$ to its inverse. We can extend this to a bijection $r_{q'}\mapsto r_{(q')^{-1}}$ since $r_{q'}$ is fully determined by the position of $i$ in $q'$. Note that $r_{q'}=(i,i-1)\cup q'$ and similarly for $r_{(q^{-1})'}$ (here taking $q'$ and $(q^{-1})'$ to be partial orders on $Q\cup\{i-1\}$).
    
    Consider the general form of $r_{q'}$. The shape of the Hasse diagram of $r_{q'}$ is that of a linear order of length $n-i$ (which is $q'$) with an `offshoot' of length $1$ at some vertex (which is $i$). We can partition the $r_{q'}$ into $n-i$ sets by the position of $i$ in $q'$ (since $a_i=-1$ we must have $(i+1,i)\in q'$ so $i$ cannot be minimal, but it can be in any other position). In effect, we are partitioning the $r_{q'}$ into their isomorphism classes; two partial orders are isomorphic if and only if their Hasse diagrams are `the same' once labels have been removed. Let $P_1,\ldots,P_{n-i}$ be the partition, where $P_j=\{q'\in L(q)\ |\ q'(i)=j+1\}$ and let $N_j=|P_j|$ for all $j$.
    
    We can also partition $\bigcup_{q'\in P_j}L(r_{q'})$ based on isomorphism class. Let $S_j^k=\{t\in\bigcup_{q'\in P_j}L(r_{q'})\ |\ t(i-1)=k+2\}$ for $j\leq k\leq n-i$. For a given $q'\in P_j$ there is precisely one $t\in L(r_{q'})$ with $t(i-1)=k+2$ for each $j\leq k\leq n-i$. Hence $|S_j^k|=N_j$ for each $j\leq k\leq n-i$. Define $S^k=\bigcup_{1\leq j\leq k}S_j^k$, this is the set of all $t\in\bigcup_{q'\in L(q)}L(r_{q'})$ with $t(i-1)=k+2$. Thus $\bigcup_{k=1}^{n-i}S^k=\bigcup_{k=1}^{n-i}\{t\in\bigcup_{q'\in L(q)}L(r_{q'})\ |\ t(i-1)=k+2\}=\bigcup_{q'\in L(q)}L(r_{q'})$ and so for $t,u\in\bigcup_{q'\in L(q)}L(r_{q'})$ refinements $p\cup t$ and $p\cup u$ are isomorphic if and only if $t,u\in S^k$ which holds if and only if $|L(p\cup t)|=|L(p\cup u)|$ and so we denote this common value by $|L(p\cup t)|_k$. Since the $P_j$ are disjoint so are $S_j^k$, thus $|S^k|=\sum_{j=1}^k|S_j^k|=\sum_{j=1}^kN_j$.
    
    With all this it follows by Lemma 3.8 that
    \begin{multline*}
        L(p)=\bigsqcup_{q'\in L(q)}L(p\cup q')=\bigsqcup_{q'\in L(q)}\bigsqcup_{t\in L(r_{q'})}L(p\cup q'\cup t)\\=\bigsqcup_{q'\in L(q),\ t\in L(r_{q'})} L(p\cup t)=\bigsqcup_{k=1}^{n-i}\bigsqcup_{t\in S^k}L(p\cup t)
    \end{multline*}
    
    and so
    \[|L(p)|=\sum_{k=1}^{n-i}\sum_{t\in S^k}|L(p\cup t)|=\sum_{k=1}^{n-i}(\sum_{j=1}^kN_j)|L(p\cup t)|_k\]
    
    Since $P_1,\ldots, P_{n-i}$ is a partition for $L(q)$, we can partition $L(q^{-1})$ similarly. If $q'\in P_j$ then $q'(i)=j+1$ and so $(q')^{-1}(i)=n-i-j+1$. Thus we can partition $L(q^{-1})$ into $P_1^{-1},\ldots,P_{n-i}^{-1}$ where $P_j^{-1}=\{(q^{-1})'\in L(q^{-1})\ |\ (q^{-1})'(i)=n-i-j+1\}$. By Lemma 3.7 $|P_j^{-1}|=|P_j|=N_j$.
    
    Continuing, define $(S_j^k)^{-1}=\{t\in\bigcup_{(q^{-1})'\in P_j^{-1}}L(r_{(q^{-1})'})\ |\ t(i-1)=k+2\}$ for $n-i-j\leq k\leq n-i$. Similar to before, for $(q^{-1})'\in P_j^{-1}$ there is precisely one $t\in L(r_{(q^{-1})'})$ with $t(i-1)=k+2$ for each $n-i-j\leq k\leq n-i$, and so $|(S_j^k)^{-1}|=N_j$. Define $(S^k)^{-1}=\bigcup_{n-i-k\leq j\leq n-i}(S_j^k)^{-1}$ taking $(S_0^{n-i})^{-1}=\emptyset$. We then have $\bigcup_{k=0}^{n-i}(S^k)^{-1}=\bigcup_{(q^{-1})'\in L(q^{-1})}L(r_{(q^{-1})'})$.
    
    The sets $(S_j^k)^{-1}$ are disjoint so $|(S^k)^{-1}|=\sum_{j=n-i-k}^{n-i}|(S_j^k)^{-1}|$ $=\sum_{j=n-i-k}^{n-i}N_j$ for $k< n-i$ and $|(S^{n-i})^{-1}|=\sum_{j=1}^{n-i}N_j$ and as before $t,u\in (S^k)^{-1}$ if and only if $p\Phi_i\cup t$ and $p\Phi_i\cup u$ are in the same isomorphism class and thus $|L(p\Phi_i\cup t)|=|L(p\Phi_i\cup u)|$. Importantly, we also have that if $t\in S^k$ and $u\in (S^k)^{-1}$ then $p\cup t$ and $p\Phi_i\cup u$ are isomorphic since $p$ and $p\Phi_i$ are equal on $\{1,\ldots,i\}$ and so $|L(p\cup t)|=|L(p\Phi_i\cup u)|$. Hence $|L(p\Phi_i\cup u)|=|L(p\cup t)|_k$ for all $u\in (S^k)^{-1}$.
    
    We now have
    \[L(p\Phi_i)=\bigsqcup_{(q^{-1})'\in L(q^{-1}),\ t\in L(r_{(q^{-1})'})}L(p\Phi_i\cup t)=\bigsqcup_{k=0}^{n-i}\bigsqcup_{t\in (S^k)^{-1}}L(p\Phi_i\cup t)\]
    
    and then 
    \begin{multline*}
        |L(p\Phi_i)|=\sum_{k=0}^{n-i}\sum_{t\in (S^k)^{-1}}|L(p\Phi_i\cup t)|\\=\sum_{k=0}^{n-i-1}(\sum_{j=n-i-k}^{n-i}N_j)|L(p\cup t)|_k+(\sum_{j=1}^{n-i}N_j)|L(p\cup t)|_{n-i}
    \end{multline*}
    
    Now, comparing the two formulae obtained, the statement of the proposition is equivalent to
    
    \[\sum_{k=1}^{n-i-1}(\sum_{j=1}^kN_j)|L(p\cup t)|_k < \sum_{k=0}^{n-i-1}(\sum_{j=n-i-k}^{n-i}N_j)|L(p\cup t)|_k\]
    
    Thus it suffices to show that $\sum_{j=1}^kN_j\leq\sum_{j=n-i-k}^{n-i}N_j$ for all $1\leq k\leq n-i-1$. As noted $q\cong Z_{n-i+1}^{-1}$ with the labelling from $i$ to $n$. Thus $i$ takes the place of $1$ in Lemma 3.9 and so $P_j\cong L(Z_{n-i+1}^{-1},j)$ and hence by that lemma we have $N_j=E(n-i,j)$. Finally note that the recurrence which defines the Entringer numbers rearranges to give $E(n,i)-E(n,i-1)=E(n-1,n-i)>0$ and thus $N_1<N_2<\ldots<N_{n-i}$. \qed
\end{pf}

The theorem now follows.

\begin{theorem}
    The canonical cycle $(1,3,5,\ldots,6,4,2)$ and its inverse have the highest multiplicity among cycles in $C(P_n)$ and are the only cycles to obtain this maximum.
\label{t:pmx}
\end{theorem}

\begin{pf}
    Let $p=(a_1,\ldots,a_{n-1})\in\mathbf{P}_n\setminus \{Z_n,Z_n^{-1}\}$. Suppose $i_1,\ldots,i_m$ are the values such that $a_{i_j}=a_{i_j-1}$ in ascending order. Then by Proposition 3.9 we have $|L(p)|<|L(p\Phi_{i_m})|<|L(p\Phi_{i_m}\Phi_{i_{m-1}})|<\ldots<|L(p\Phi_{i_m}\ldots\Phi_{i_1})|$ and either $p\Phi_{i_m}\ldots\Phi_{i_1}\cong Z_n$ or $p\Phi_{i_m}\ldots\Phi_{i_1}\cong Z_n^{-1}$. Since $p$ was arbitrary the result follows. \qed
\end{pf}

Notice that in the proof of Proposition 3.10 the assumption that $i$ is maximal is not used until the final paragraph. With this in mind we make the following conjecture.

\begin{conjecture}
    Let $p=(a_1,\ldots,a_{n-1})\in\mathbf{P}_n$ and let $i\geq 2$ be such that $a_{i-1}=a_i$. Then $|L(p)|<|L(p\Phi_i)|$.
\end{conjecture}

Since right-side inversion is an involution, proving this would mean that given any $p\in\mathbf{P}_n$ and any $2\leq i\leq n-1$ we would know which of $p$ and $p\Phi_i$ has more linear extensions. Identifying elements which are inverse (they have the same number of linear extensions) we can obtain a bounded lattice $L$ on $\mathbf{P}_n$ by defining $(p, p\Phi_i)\in L$ if the above condition holds and taking the transitive closure, and so if the conjecture is true then $L$ has the property that $|L(p)|<|L(q)|$ if $(p,q)\in L$. This lattice $L$ is isomorphic to the power set of a set with $n-2$ elements ordered by inclusion.

\subsubsection{The least frequent cycles}

\begin{prop}
For the path with vertices numbered from $1$ to $n$ in order, the least
frequent cycles are $(1,2,\ldots,n)$ and its inverse, which are realised
just once.
\end{prop}

\begin{pf}
The count follows from Corollary~\ref{c:uniq_cyc}; the proof of that Corollary
shows that to realise these cycles we must take the edges of the path in order
from one end to the other. Now
\[(1,2)(2,3)\cdots(n-1,n)=(n,n-1,\ldots,2,1),\]
and similarly for the reverse order. \qed
\end{pf}

The second smallest frequency appears to be $n-2$, realised by four cycles if
$n>4$. To prove this, we first require a couple of lemmas.

\begin{lemma}
    Let $p=(a_1,\ldots,a_{n-1})\in\mathbf{P}_n$ and let $1\leq i\leq n-2$ be such that $a_{i+j}=1$ for all $j\geq 0$ or $a_{i+j}=-1$ for all $j\geq 0$. Then $|L(p)|<|L(p\Phi_{i+j})|$ for all $j\geq 1$.
\end{lemma}

\begin{pf}
    Consider the induced suborder $q$ of $p$ on $Q=\{i,\ldots,n\}$. This is a linear order by assumption and thus has precisely one linear extension. On the other hand, given $j\geq 1$ the induced suborder $q_j$ of $p\Phi_{i+j}$ on $Q$ is not a linear order and thus has at least two linear extensions. In particular, there is a linear extension $r$ of $q_j$ such that $p\Phi_{i+j}\cup r$ is isomorphic to $p\cup q=p$. If $a_i=1$ (resp. $a_i=-1$) then any $r\in L(q_j)$ such that $r(i)=1$ (resp. $r(i)=n-i+1$) is appropriate. By Lemma 3.8 we have \[L(p\Phi_{i+j})=\bigsqcup_{q_j'\in L(q_j)}L(p\Phi_{i+j}\cup q_j')\supset L(p\Phi_{i+j}\cup r)\] and thus $|L(p)|=|L(p\Phi_{i+j}\cup r)|<|L(p\Phi_{i+j})|$. \qed
\end{pf}

The next lemma concerns path-type partial orders of the form $(1,\ldots,1,$ $-1,\ldots,-1)\in\mathbf{P}_n$. If $(a_1,\ldots,a_{n-1})$ is such a partial order we denote it by $h_k$, where $k$ is the unique number such that $a_{k-1}=1$ and $a_k=-1$. Then the vertex of $h_k$ labelled $k$ is the global maximal element of $h_k$. We denote $\mathbf{H}_n=\{h_k\in\mathbf{P}_n\ |\ k=2,\ldots,n-1\}$ and refer to this set's elements as \emph{hills}; for a hill $h_k$ we refer to the vertex labelled $k$ as its \emph{peak}.

\begin{lemma}
    Consider $\mathbf{H}_n$ for some $n$. We have for all $k$ \[|L(h_k)|=
    {n-1\choose n-k}.\]
\end{lemma}

\begin{pf}
    Suppose $i>\frac{n}{2}$, the complementary case is symmetrical and follows by the same argument. We count the linear extensions `manually'; let $h\in L(h_k)$, this will always be the linear extension we are currently considering. We take $h(1)<h(2)<\ldots<h(k)$ to be fixed reference points and define the extension we are considering by placing the remaining points $h(k+1)>\ldots>h(n)$ in between these fixed points; each choice uniquely defines a linear extension and every linear extension can be defined this way. However these choices are not free. We have $h(i)>h(i+1)$ for $i=k+1,\ldots,n-1$ so choosing where $h(i)$ goes restricts the choices for $h(i+1)$ and thus we choose where $h(n)$ goes first, then $h(n-1)$ and so on. Further, $h(k)>h(k+1)$ so there are $k$ places to choose from. We denote by $G_i$ the space between $h(i-1)$ and $h(i)$ for $i=2,\ldots k$ and by $G_1$ the space below $h(1)$.
    
    We start with $h(k+1),\ldots,h(n)\in G_1$, this is one linear extension. Keeping $h(k+2),\ldots,h(n)\in G_1$ and placing $h(k+1)\in G_2$ is another linear extension. Similarly $h(k+2),\ldots,h(n)\in G_1$ and $h(k+1)\in G_i$ for $i=3,\ldots,k$ each gives a linear extension and we have found $\sum_{i=1}^k1$ linear extensions so far. Next, we start with $h(k+1),h(k+2)\in G_2$ and $h(k+3),\ldots,k(n)\in G_1$ and we find a linear extension $h(k+1)\in G_i$, $h(k+2)\in G_2$, $h(k+3),\ldots,k(n)\in G_1$ for each $i=3,\ldots,k$, giving us another $\sum_{i=1}^{k-1}1$ linear extensions. Moving $h(k+2)$ in the same way as we have with $h(k+1)$ and then moving $h(k+1)$ within the thus restricted zone, we find a linear extension $h(k+1)\in G_i$, $h(k+2)\in G_j$, $h(k+3),\ldots,h(n)\in G_1$ for each $i,j=3,\ldots,k$ such that $i\geq j$. Choosing $j$ first then finding each $i$ within the possible range, for such a $j$ we have $\sum_{i=1}^{k-j+1}1$ linear extensions. Hence, in total, we have now found $\sum_{i_1=1}^k\sum_{i_0=1}^{i_1}1$ linear extensions.
    
    Continuing inductively, we find
    that \[|L(h_k)|=\sum_{i_{n-k-1}=1}^k\sum_{i_{n-k-2}=1}^{i_{n-k-1}}\ldots\sum_{i_1=1}^{i_2}\sum_{i_0=1}^{i_1}1.\] Clearly $\sum_{i_0=1}^{i_1}1=i_1={i_1\choose 1}$ and it is well known that \[\sum_{i_1=1}^{i_2}i_1=\frac{i_2(i_2+1)}{2}={i_2+1\choose 2}.\] It follows from the hockey stick identity \cite{jones} that for $k\geq 0$ \[\sum_{i=1}^m{i+k-1\choose k}={m+k\choose k+1}\] and thus by induction on $m$ we have \[\sum_{i_m=1}^{i_{m+1}}\sum_{i_{m-1}=1}^{i_m}\ldots\sum_{i_1=1}^{i_2}i_1={i_{m+1}+m\choose m+1}.\] Finally, we conclude \[|L(h_k)|=\sum_{i_{n-k-1}=1}^k{i_{n-k-1}+n-k-2\choose n-k-1}={n-1\choose n-k}.\] \qed
\end{pf}

\begin{prop}
    Consider the path with $n+1$ vertices labelled consecutively from $1$ to $n+1$. If $n\geq4$ then $(1,2,\ldots,n-1,n+1,n)$ and $(1,3,4,\ldots,n,n+1,2)$ and their inverses are the second least frequent cycles, appearing with multiplicity $n-1$.
\end{prop}

\begin{pf}
    Notice that for $n=1,2$ there are only one and two distinct cycles respectively. For $n=3$ there are four distinct cycles; two least frequent cycles and two most frequent cycles. So the first interesting case here is $n=4$. 
    
    Let $e_1,\ldots,e_n$ be the edges of the path such that $e_i=(i,i+1)$. Then we have \[\sigma=e_1e_2\ldots e_{n-2}e_ne_{n-1}=(1,2,\ldots,n-1,n+1,n)^{-1}\] and \[\tau=e_2e_1e_3\ldots e_{n-1}e_n=(1,3,4,\ldots,n,n+1,2)^{-1}\]
    
    Thus we can use these two orderings of the edges to find the partial orders corresponding to these cycles. We find that the local suborders of $\sigma$ are $\sigma|_i=e_{i-1}e_i$ for $2\leq i\leq n-1$ and $\sigma|_{n}=e_{n}e_{n-1}$. Similarly, we have $\tau|_2=e_2e_1$ and $\tau|_i=e_{i-1}e_i$ for $3\leq i\leq n$. Thus the corresponding partial orders are $p_\sigma=(1,\ldots,1,-1),p_\tau=(-1,1,\ldots,1)\in\mathbf{P}_n$ whose Hasse diagrams are shown in Figure~\ref{f:2lf}.
    
    \begin{figure}[htbp]
        \begin{center}
            \setlength{\unitlength}{1mm}
            \begin{picture}(50,28)
                \multiput(2,2)(3,6){3}{\circle*{1}}
                \multiput(11,20.5)(3,6){2}{\circle*{1}}
                \put(17,20.5){\circle*{1}}
                \put(2,2){\line(1,2){6}}
                \put(11,20.5){\line(1,2){3}}
                \put(8.5,15){\rotatebox{65}{\dots}}
                \put(17,20.5){\line(-1,2){3}}
                \put(1,-5){$p_{\sigma}=h_{n-1}$}
                
                \multiput(39,2)(3,6){2}{\circle*{1}}
                \multiput(45,14.5)(3,6){3}{\circle*{1}}
                \put(36,8){\circle*{1}}
                \put(36,8){\line(1,-2){3}}
                \put(39,2){\line(1,2){3}}
                \put(45,14.5){\line(1,2){6}}
                \put(42.5,9){\rotatebox{65}{\dots}}
                \put(35,-5){$p_{\tau}=h_2^{-1}$}
            \end{picture}
        \end{center}
        \caption{\label{f:2lf}The partial orders of the second least frequent cycles for paths}
    \end{figure}
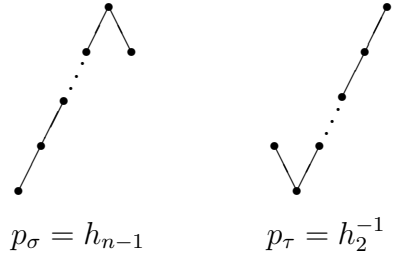
    
    Notice that $p_\sigma=h_{n-1}$ and $p_\tau=h_2^{-1}$; thus by Lemma 3.15 we have $|L(p_\sigma)|={n-1\choose 1}={n-1\choose n-2}=|L(p_\tau)|=n-1$ and so this is the multiplicity of the cycles. To see that these are the second least frequent cycles, let $p=(a_1,\ldots,a_{n-1})\in\mathbf{P}_n$ such that $p\neq(1,\ldots,1),(-1,\ldots,-1),h_2^{\pm1},h_{n-1}^{\pm1}$ and let $i_1,\ldots,i_m$ be the numbers such that $a_{i_j-1}=-a_{i_j}$ listed in ascending order. Then $p\cong h_{i_1}^{a_1}\Phi_{i_2}\ldots\Phi_{i_m}$ and by Lemma 3.14 we have $|L(h_{i_1}^{a_1})|<|L(h_{i_1}^{a_1}\Phi_{i_2})|<\ldots<|L(h_{i_1}^{a_1}\Phi_{i_2}\ldots\Phi_{i_m})|=|L(p)|$. Finally by Lemma 3.15 we have $|L(h_{n-1})|=|L(h_2)|<|L(h_k)|$ for any $2<k<n-1$ and we conclude $|L(h_2)|=|L(h_{n-1
    })|<|L(p)|$. \qed
\end{pf}

\paragraph{Problem} For the path on $n$ vertices, describe the cycles with the second largest number of realisations, and
calculate this number.

\subsection{Forked paths}

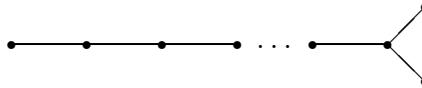
\begin{figure}[htbp]
\begin{center}
\setlength{\unitlength}{1mm}
\begin{picture}(60,10)
\multiput(5,5)(10,0){6}{\circle*{1}}
\put(5,5){\line(1,0){30}}
\put(37.5,4.5){\dots}
\put(45,5){\line(1,0){10}}
\multiput(60,0)(0,10){2}{\circle*{1}}
\put(55,5){\line(1,1){5}}
\put(55,5){\line(1,-1){5}}
\end{picture}
\end{center}
\caption{\label{f:fork}A forked path}
\end{figure}

We have done limited investigation of other trees. One natural candidate is
the ``forked path'' or Coxeter diagram of type $D_n$ for $n\ge3$
(Figure~\ref{f:fork}).

Empirically the greatest frequency of cycles obtained from this tree is given
by Sequence A034428 in the OEIS, counting \emph{almost up-down permutations}
(which begin with two increases and then alternate); they are given by the
formula $A_n=nE_{n-1} - E_{n}$, where $(E_n)$ are the Euler numbers.
The sequence begins
\[1,1,3,9,35,155,791,4529,28839,\ldots\]
and has generating function $1-(1-x)(\tan(x)+\sec(x))$. The analogous
almost down-up permutations and their reverses also give the same number of
cycles.

First, we give a characterisation of $\mathbf{D}_n=P(D_{n+1})$ similar to that of Proposition 3.6. We assume the leaves at the `fork' are labelled $1$ and $2$ and the rest of the vertices labelled $3,\ldots,n+1$ consecutively, while the edges are labelled $e_i$ so that $e_1=(1,3)$ and $e_i=(i,i+1)$ for $i=2,\ldots,n$. For convenience we will often refer to $e_i$ simply as $i$ when the referent is clear; in particular when referring to labels on the vertices of Hasse diagrams.

\begin{prop}
    The set $\mathbf{D}_n$ is precisely the set of partial orders whose Hasse diagram is a path-type partial order with $n-2$ vertices with two vertices attached to the leftmost vertex making a linear suborder on three vertices.
\end{prop}

\begin{pf}
    Take some $\sigma\in O(D_{n+1})$ and consider $p_\sigma\in\mathbf{D}_n$. As with path-type partial orders we find that the vertices $1,2$ and $n+1$ are leaves while $\sigma|_i=e_{i-1}e_i$ or $e_ie_{i-1}$ for $i=4,\dots,n$ so the induced suborder of $p_\sigma$ on $\{3,\ldots,n\}$ is a path-type partial order. The final part to consider is $\sigma|_3$, where $e_1,e_2,e_3$ appear in some order; $e_3$ is the only one of these to appear in any other cover relation, so the Hasse diagram of $p_\sigma$ is found by attaching the vertices $1$ and $2$ to $3$ so that they form a linear suborder.
    
    Now, by Theorem~\ref{t:t2c} we know that $|\mathbf{D}_n|=|C(D_{n+1})|=6\times 2^{n-3}$. We show that this is the same as the number of Hasse diagrams of the stated type. Pick a path-type Hasse diagram on $n-2$ vertices labelled consecutively from $3$ to $n$ (there are $2^{n-3}$ of them) and a linear order on $1,2,3$ (there are $3!=6$ of them) and attach the linear order to the leftmost vertex of the path by identifying the vertices labelled $3$. This fully determines a Hasse diagram of the stated type, and each such diagram can be reached in this way. \qed
\end{pf}

We may refer to these as partial orders of \emph{fork-type}; given a fork-type partial order we refer to the linear suborder on $\{1,2,3\}$ as the \emph{fork part} and the induced suborder on $\{3,\ldots,n\}$ as the \emph{path part}. As such given a fork-type partial order $d\in\mathbf{D}_n$ we will often denote it as $(a_1,a_2,a_3)p$ where $(a_1,a_2,a_3)$ tells us that the fork part has cover relations $(a_1,a_2),(a_2,a_3)$ and $p\in\mathbf{P}_{n-2}$ is the path part with relabelling $i\mapsto i+2$. Given this, we note that Proposition 3.10 can be extended to $\mathbf{D}_n$ as long as right-side inversion is applied at a point in the path part, and we will need this fact in the proofs that follow.

We can now use this characterisation to find the most frequent cycles and their multiplicity as we did with paths. For this purpose we first find an explicit formula for the number of linear extensions of a particular fork-type partial order.

\begin{lemma}
    Consider $d=(a_1,3,a_2)Z_{n}^{a}\in\mathbf{D}_{n+2}$ for $a=\pm 1$ and some $n\in\mathbb{N}$. Then \[|L(d)|=\sum_{i=1}^{n-1}(n-i)(i+1)E(n-1,i).\]
\end{lemma}

\begin{pf}
    We assume $d=(1,3,2)Z_{n}^{-1}$; the other cases follow by symmetry and inverting. The formula follows by considering refinements of $d$ found by taking a linear extension of the path part and summing the number of linear extensions for each refinement, since by Lemma 3.8 \[|L(d)|=\sum_{z\in L(Z_n^{-1})}|L(d\cup z)|\] recalling $Z_n^{-1}$ is relabelled by $i\mapsto i+2$. 
    
    The Hasse diagram of a refinement of this type is a linear order on three points `joined' with a linear order on $n$ points by identifying the vertices labelled by $3$. Thus the isomorphism class of a refinement $d\cup z$ is fully determined by the place of the vertex labelled $3$ in the linear extension $z$. As such we partition the refinements $\bigcup_{z\in L(Z_n^{-1})}\{d\cup z\}$ into sets $D_i=\{d\cup z\ |\ z\in L(Z_n^{-1},i)\}$ for $i=1,\ldots, n-1$ and any two refinements in a given set are isomorphic.
    
    Consider $r\in D_i$ for some $i$. The relationship of the vertex labelled $3$ with each other vertex is already fully determined, so when finding a linear extension of $r$ this cannot be altered. Thus the linear extensions of $r$ can be counted by considering it as a hill $b=h_2\in\mathbf{H}_{i+2}$ and an inverted hill $t=h_2^{-1}\in\mathbf{H}_{n-i+1}^{-1}$ joined by identifying their peaks, which is the vertex labelled $3$. This is because if $a\in[n+2]$ is such that $(a,3)\in r$ then a linear extension $r'$ cannot have $(3,a)\in r'$, and similarly if $a$ is such that $(3,a)\in r$. Hence any linear extension of $r$ can be found by first taking a linear extension $b'$ of $b$ and then of $r\cup b'$, which has the same number of linear extensions as $t$. Thus we have $|L(r)|={i+1\choose i}{n-i\choose n-i-1}=(i+1)(n-i)$ by Lemma 3.15. Finally, continuing the equation above we have \[|L(d)|=\sum_{i=1}^{n-1}\sum_{z\in D_i}(n-i)(i+1)=\sum_{i=1}^{n-1}(n-i)(i+1)E(n-1,i)\] by Lemma 3.9. \qed
\end{pf}

\begin{theorem}
    The cycles $(1,2,4,6,\ldots,7,5,3),(1,2,3,5,7,\ldots,6,4)$ and their inverses are the most frequent cycles among cycles in $C(D_{n+1})$ having multiplicity $A_n$, and are the only such to attain this multiplicity.
\end{theorem}

\begin{pf}
    Notice that \[\sigma=e_1e_2e_4e_6\ldots e_5e_3=(1,2,4,6,\ldots,7,5,3)\] and \[\tau=e_2e_1e_4e_6\ldots e_5e_3=(1,2,3,5,7,\ldots,6,4)^{-1}\] so we can find the multiplicity of these cycles by considering the partial orders $p_\sigma,p_\tau\in\mathbf{D}_n$. Calculating each of these we see that they are isomorphic; their isomorphism class is shown in Figure~\ref{f:az}.
    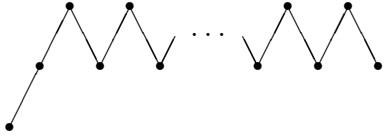
\begin{figure}[htbp]
        \begin{center}
            \setlength{\unitlength}{1mm}
            \begin{picture}(60,17)
                \multiput(5,-2)(4,8){3}{\circle*{1}}
                \multiput(17,6)(4,8){2}{\circle*{1}}
                \put(25,6){\circle*{1}}
                \put(5,-2){\line(1,2){8}}
                \put(13,14){\line(1,-2){4}}
                \put(17,6){\line(1,2){4}}
                \put(21,14){\line(1,-2){4}}
                \put(25,6){\line(1,2){2}}
                \put(29,10){\dots}
                \put(36,10){\line(1,-2){2}}
                \multiput(42,14)(8,0){2}{\circle*{1}}
                \multiput(38,6)(8,0){3}{\circle*{1}}
                \put(38,6){\line(1,2){4}}
                \put(42,14){\line(1,-2){4}}
                \put(46,6){\line(1,2){4}}
                \put(50,14){\line(1,-2){4}}
            \end{picture}
        \end{center}
    \caption{\label{f:az}The almost zig-zag partial order for even $n$}
    \end{figure}
    Using an analogous argument to that used in Lemma 3.9 we can see that the number of linear extensions of this partial order is the same as the number of almost up-down permutations on $n$ points, which is $A_n$ by definition.
    
    We may partition $\mathbf{D}_n$ into six sets defined by \[\mathbf{D}_n^{(a_1,a_2,a_3)}=\{d\in\mathbf{D}_n\ |\ d=(a_1,a_2,a_3)p\ \text{for some} \ p\in\mathbf{P}_{n-2}\}.\] We simplify this by noting that $\mathbf{D}_n^{(a_1,a_2,a_3)}$ is precisely the set of inverses of $\mathbf{D}_n^{(a_3,a_2,a_1)}$ and so we need only consider three sets of the partition - $\mathbf{D}_n^{(1,2,3)},\mathbf{D}_n^{(1,3,2)}$ and $\mathbf{D}_n^{(3,1,2)}$.
    
    We begin by noting that the partial orders with the most linear extensions in each of these sets are the ones whose path part is zig-zag. Indeed, given $a_1,a_2,a_3$ from one of the three above and $p\in\mathbf{D}_n^{(a_1,a_2,a_3)}$ we may repeatedly apply right-side inversion in the same manner as in Theorem~\ref{t:pmx} to see that it has strictly fewer linear extensions than an appropriate such partial order. The path part can be isomorphic to either $Z_{n-2}$ or $Z_{n-2}^{-1}$ and so there are six such partial orders in $\mathbf{D}_n^{(1,2,3)}\cup\mathbf{D}_n^{(1,3,2)}\cup\mathbf{D}_n^{(3,1,2)}$ which are $(a_1,a_2,a_3)Z_{n-2},(a_1,a_2,a_3)Z_{n-2}^{-1}\in\mathbf{D}_n^{(a_1,a_2,a_3)}$. Notice that $p_\sigma=(1,2,3)Z_{n-2}^{-1}$ and $p_\tau=(3,1,2)Z_{n-2}$.
    
    We can see that $p_\sigma=(1,2,3)Z_{n-2}\Phi_3$ and $p_\tau=(3,1,2)Z_{n-2}^{-1}\Phi_3$ and thus $|L(p_\sigma)|>|L((1,2,3)Z_{n-2})|$ and $|L(p_\tau)|>|L((3,1,2)Z_{n-2}^{-1})|$. All that is left to show is $|L((1,3,2)Z_{n-2})|,|L((1,3,2)Z_{n-2}^{-1})|<|L(p_\sigma)|=|L(p_\tau)|$ for $n\geq 4$. By Lemma 3.18 we have \[|L((1,3,2)Z_{n-2})|=|L((1,3,2)Z_{n-2}^{-1})|=\sum_{i=1}^{n-3}(n-i-2)(i+1)E(n-3,i)\] and using $E(n-1,n-i)=E(n,i)-E(n,i-1)$ we can rearrange this as 
    \begin{multline*}
        =\sum_{i=1}^{k}(n-2(i+2))[E(n-2,n-i-3)-E(n-2,i+1)] \\
        + 2(n-3)E(n-2,n-3) - (n-4)E(n-2,1)
    \end{multline*}
    where $k=\frac{n}{2}-3$ if $n$ is even and $k=\frac{n+1}{2}-3$ if $n$ is odd. Taking the difference between $A_n=nE_{n-1}-E_n$ and this by using $E(n,i)=\sum_{k=n-i}^{n-1}E(n-1,k)$ and $E_n=E(n,n)$ we find
    \begin{multline*}
        \sum_{i=1}^{k+1}(2n-3i-3)E(n-2,i) + \sum_{i=1}^{k}(3i+6-n)E(n-2,n-i-3) \\
        + (8-2n)E(n-2,n-3) + E(n-2,n-2) \\
        + (1 - (n\mod2) )\frac{n}{2}E(n-2,\frac{n}{2}-1)
    \end{multline*}
    and thus it suffices to show that this sequence is positive for all $n\geq 4$. By examination we can see that this is equal to 
    \[\sum_{i=1}^{n-4}(2n-3i-3)E(n-2,i)-(2n-8)E(n-2,n-3)+E(n-2,n-2)\] for all $n$. Applying $E(n-2,i)=\sum_{k=n-i-2}^{n-3}E(n-3,k)$ and rearranging to express as a linear combination of all the $(n-3)$th Entringer numbers we obtain
    \begin{multline*}
        \sum_{i=1}^{n-3}\Bigg(\sum_{k=n-i-2}^{n-4}(2n-3k-3)-2n+9\Bigg)E(n-3,i) \\
        =\sum_{i=1}^{n-3}\frac{3}{2}\big(i^2+(3-\frac{2}{3}n)i+(2-\frac{2}{3}n)\big)E(n-3,i)
    \end{multline*}
    Let $C_i=\frac{3}{2}(i^2+(3-\frac{2}{3}n)i+(2-\frac{2}{3}n))$ for $i=0,\ldots,n-3$. Then we calculate that $\sum_{i=1}^{n-3}C_i=n-3$ and by considering $C_i$ as a polynomial on $i$ we can further deduce that $C_i>0$ precisely when $i>\frac{2}{3}n-2$. This tells us that
    there are more positive terms than negative (treating $C_iE(n-3,i)$ as $|C_i|$ instances of the term $E(n-3,i)$) and the smallest positive term is greater than the largest negative term (since $E(n-3,j)<E(n-3,k)$ for $j<k$). Thus the sequence is positive as required. \qed
\end{pf}

\section{The inverse problem}

Now we turn to the question: how many labelled trees yield a given cycle $c$
when the edges are multiplied together? Since the cycles are all conjugate in
the symmetric group, the answer is independent of the chosen cycle; so where
necessary we can assume that $c=(1,2,\ldots,n)$.

By Theorem~\ref{t:t2c}, if the vertices of $T$ have valencies $d_1,\ldots,d_n$,
then the number of cycles which can be obtained from $T$ is
\[C(T)=\prod_{i=1}^nd_i!.\]
The symmetric group acts on the set of $n^{n-2}$ trees. This action is not
transitive (for $n>2$); the orbits are the isomorphism types of trees. The
stabiliser of a tree $T$ in this action is its automorphism group $\Aut(T)$, so
the number of trees in the isomorphism class of $T$ is $n!/|\Aut(T)|$, by the
Orbit-Stabiliser Theorem.

If each cycle has $D(T)$ realising trees in the isomorphism class of $T$ (the
orbit of $S_n$ containing $T$), then we have
\[(n-1)!D(T)=(n!/|\Aut(T)|)C(T),\]
so that $D(T)=nC(T)/|\Aut(T)|$. Summing over the isomorphism types $T$ gives
the answer to our question.

For example, for $n=5$, there are three isomorphism classes of trees, with
automorphism groups of orders $24$, $2$, and $2$, with vertex valencies other
than $1$ being respectively $(4)$, $(3,2)$, and $(2,2,2)$; so the number is
\[5\cdot24/24+5\cdot12/2+5\cdot8/2=55.\]

In order to find a more explicit formula, we begin with a couple of definitions.

\paragraph{Definition}
Let $T$ be a tree on the vertex set $\{1,\ldots,n\}$, and $c$ a cyclic
permutation of the vertices. The \emph{diagram} of the pair $(T,c)$ consists
of $n$ points $a_1,\ldots,a_n$ in the order around a circle given by the cycle
$c$, with a line segment from $a_i$ to $a_j$ whenever $\{i,j\}$ is an edge of
$T$. The diagram is \emph{crossing-free} if no pair of these line segments
intersect (except at their ends).

\paragraph{Definition}
The pair $(T,c)$ (as above) is \emph{realisable} if the product of the
transpositions corresponding to edges of $T$ (in some order) is $c$.

\begin{theorem}
The pair $(T,c)$ is realisable if and only if its diagram is crossing-free.
\label{t:c2t}
\end{theorem}

To prove the forward implication, we require the following lemma.

\begin{lemma}
Suppose that $e=\{i,j\}$ is an edge of $T$; let $T_i$ and $T_j$ be the two
trees obtained by deleting $e$, with $i\in T_i$ and $j\in T_j$. Suppose that
the product of the edges of $T$ in some order maps $i$ to $j$. Then we can
rearrange the product as $w_jew_i$, where $w_i$ is a product of the edges of
$T_i$ in some order, and so is a cycle on the vertices of $T_i$; and similarly
for $w_j$ and $T_j$.
\end{lemma}

\begin{pf} Let $u_1eu_2=(a_1,\ldots,a_r,i,j,b_1,\ldots,b_s)$. 
Then $i$ is fixed by $u_1$, so the edges of $i$ in $T_i$ are on the
right of $e$; and similarly, edges of $j$ in $T_j$ are on the left of $e$.
Edges in $T_j$  which occur to the right of $e$ in the product can thus be
moved, one by one, to the left (they commute with edges of $T_i$ and with
$(i,j)$ by the preceding remark), and vice versa. So we can end up with the
edges of $T_i$ on the right and those of $T_j$ on the left. Now call this
product $w_jew_i$. The claims of the lemma now follow. \qed
\end{pf}

Now we turn to the proof of the ``only if'' direction in the theorem. This is
a proof by induction. The result is easily verified for small values of $n$,
starting the induction. So we assume that $(T,c)$ has a crossing-free diagram,
and that all smaller pairs with crossing-free diagrams are realisable.

Let $e=\{i,j\}$ be an edge. Because the diagram is crossing-free, it falls into
two parts, one on each side of $e$.

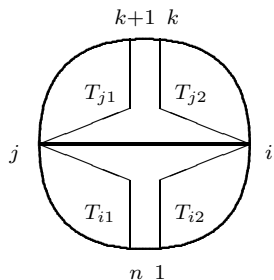
\begin{figure}[htbp]
\begin{center}
\setlength{\unitlength}{2mm}
\begin{picture}(30,30)
\thicklines
\closecurve(15,8,22,15,15,22,8,15)
\put(8,15){\line(1,0){14}}
\put(14,6){$\scriptstyle{n\,\,\,1}$}
\put(6,14){$\scriptstyle{j}$}
\put(23,14){$\scriptstyle{i}$}
\put(13,23){$\scriptstyle{k+1\,\,\,k}$}
\thinlines
\multiput(14,8)(2,0){2}{\line(0,1){4.6}}
\put(8,15){\line(5,-2){6}}
\put(22,15){\line(-5,-2){6}}
\multiput(14,22)(2,0){2}{\line(0,-1){4.6}}
\put(8,15){\line(5,2){6}}
\put(22,15){\line(-5,2){6}}
\put(11,10){$\scriptstyle{T_{i1}}$}
\put(17,10){$\scriptstyle{T_{i2}}$}
\put(11,18){$\scriptstyle{T_{j1}}$}
\put(17,18){$\scriptstyle{T_{j2}}$}
\end{picture}
\end{center}
\caption{\label{f:ind}The inductive step}
\end{figure}

Now $e$ together with the edges on one side is a non-crossing diagram, with $e=\{i,j\}$ on the boundary. By the induction hypothesis, the edges
on each side can be ordered so that the product of those on one side is
$(j,\ldots,i)$ and those on the other is $(i,\ldots,j)$ (following the order 
in $c$). By the lemma, these products can be written as $w_1ew_2$ and
$w_3ew_4$, where $w_1$ and $w_2$ are products of the edges in the two parts
on the first side left by removing $e$, and similarly on the right. Thus,
in Figure~\ref{f:ind}, $w_1$ corresponds to the edges in $T_{i2}$, and so on.

We can choose the numbering
so that $w_1ew_2=(1,2,\ldots,i,j,j+1,\ldots,n)$, where $1,\ldots,i$ and
$j,\ldots,n$ are the vertices of the two trees on the left remaining when $e$
is removed; and we have $w_1=(j,j+1,\ldots,n)$ and $w_2=(1,2,\ldots,i)$.

Similarly, there exists $k$ such that $(k+1,\ldots,j,i,i+1,\ldots,k)=w_3ew_4$,
where $w_3=(i,i+1,\ldots,k)$ and $w_4=(k+1,\ldots,j)$ correspond to the two
trees on the other side of $e$ obtained by removing $e$.

Now calculation shows that
\[w_3w_1ew_2w_4=(1,2,\ldots,i,i+1,\ldots,k,k+1,\ldots,j,j+1,\ldots,n),\]
and we are done.

\medskip

We now prove the reverse implication in the theorem. Suppose that the diagram
of the pair $(T,c)$ is not crossing-free. Then there exists two edges which intersect in the diagram, call them $e_i=\{i_1,i_2\}$ and $e_j=\{j_1,j_2\}$ where $i_1<j_1<i_2<j_2$ on the circle.

Since $T$ is a tree there is precisely one path between $e_i$ and $e_j$ and since they are non-adjacent this path has at least one edge, so choose one and call it $e_k$. If we assume that $(T,c)$ is realisable then there must exist some ordering of the edges of $T$ whose traversal (beginning at $a_1$) will witness landing on $i_1,j_1,i_2,j_2$ in that order (possibly with others in between). However, this would require landing on a vertex of $e_i$ followed by a vertex of $e_j$, the other vertex of $e_i$ and then the other vertex of $e_j$ and thus $e_k$ must appear in this traversal at least three times, contradicting Lemma 2.1. Hence $(T,c)$ is not realisable as required. \qed

\begin{cor}
The number of trees on $n$ vertices whose edges can be ordered so as to
realise a given $n$-cycle $c$ is ${3n-3\choose n-1}/(2n-1)$.
\end{cor}

\begin{pf}
See Noy~\cite{noy} for a proof that the number of crossing-free diagrams is
as claimed; see also ~\cite{dp,hp}.\qed
\end{pf}

\paragraph{Remark} These numbers are sometimes called \emph{generalised Catalan
numbers}, or \emph{Fuss--Catalan numbers}. They occur as sequence A001764 in the
On-line Encyclopedia of Integer Sequences~\cite{oeis}. The sequence begins
\[1, 1, 3, 12, 55, 273, 1428, 7752, 43263, 246675, 1430715, \ldots\]

\end{document}